\documentclass[10pt, english]{amsart}

\usepackage{savesym} 
\usepackage{fourier}
\usepackage{mathtools}
\usepackage{stmaryrd} 
\usepackage{enumerate}
\usepackage{comment}
\usepackage{todonotes}

\usepackage{graphicx}
\usepackage{tikz}
 
\newcommand{\graficosdaprova}{
\begin{tikzpicture}[
    x=0.84cm, y=1cm, >=stealth,
    declare function={f(\x)=0.9*ln(1+\x);
                      la(\x)=0.8*\x+1.5;},
    wall/.style={black, dashed, thin},
    target/.style={black, densely dotted, thick},
    bpath/.style={black, very thick, line join=round},
    apath/.style={black, thin},
    stage/.style={black!45, dash pattern=on 2pt off 2pt, very thin}]
\def\xA{1}\def\xH{2.909}\def\xB{3.5}\def\xC{5.2}\def\xD{8.9}\def\xK{11.926}\def\xE{12.6}\def\xmax{12.6}
\fill[black!10]  (\xB,-2.7) rectangle (\xC,9.9);
\fill[black!10]  (\xD,-2.7) rectangle (\xE,9.9);
\foreach \x/\lab in {\xA/{n_{2t}},\xB/{n_{2t+1}},\xC/{n_{2t+2}},\xD/{n_{2t+3}},\xE/{n_{2t+4}}}
  {\draw[stage] (\x,-2.7) -- (\x,9.9);
   \node[below,font=\footnotesize] at (\x,-2.7) {$\lab$};}
\foreach \a/\b/\lab in {\xA/\xB/stationary,\xB/\xC/growth,\xC/\xD/stationary,\xD/\xE/growth}
  \node[above,font=\footnotesize\itshape] at ({(\a+\b)/2},9.9) {\lab};
\begin{scope}[yshift=4.6cm, y=0.35cm]
  \draw[->] (0,0) -- (\xmax+0.4,0) node[right,font=\footnotesize] {$\log n$};
  \draw[->] (0,0) -- (0,15);
  \draw[wall]  plot[domain=0:\xmax,samples=60] (\x,{la(\x)+f(\x)})
       node[right,font=\footnotesize,black] {$a(n)\log n$};
  \draw[wall]  plot[domain=0:\xmax,samples=60] (\x,{la(\x)-f(\x)})
       node[right,font=\footnotesize,black] {$a(n)/\log n$};
  \draw[apath] (0,{la(0)}) -- (\xmax,{la(\xmax)}) node[right,font=\footnotesize] {$a(n)$};
  \draw[bpath] (0,{la(0)+0.3}) -- (\xA,{la(\xA)+0.3}) -- (\xH,{la(\xH)-f(\xH)})
      -- plot[domain=\xH:\xB,samples=10] (\x,{la(\x)-f(\x)})
      -- (\xC,6.686) -- (\xD,6.686) -- (\xK,{la(\xK)+f(\xK)})
      -- plot[domain=\xK:\xE,samples=10] (\x,{la(\x)+f(\x)});
  \draw[apath] (0.4,12.6) -- (1.0,12.6) node[right,font=\footnotesize,black] {$a(n)$};
  \draw[bpath] (0.4,11.3) -- (1.0,11.3) node[right,font=\footnotesize,black] {$b(n)$};
  \node[font=\footnotesize,anchor=north west] at (0.1,14.8) {log scale};
\end{scope}
\begin{scope}
  \draw[->] (0,0) -- (\xmax+0.4,0);
  \draw[->] (0,-2.6) -- (0,2.7) node[above,font=\footnotesize] {$\log\frac{b(n)}{a(n)}$};
  \draw[wall] plot[domain=0:\xmax,samples=60] (\x,{f(\x)})
       node[right,font=\footnotesize,black] {$+\log\log n$};
  \draw[wall] plot[domain=0:\xmax,samples=60] (\x,{-f(\x)})
       node[right,font=\footnotesize,black] {$-\log\log n$};
  \draw[target] (\xA,0.9) -- (\xC,0.9)  (\xA,-0.9) -- (\xC,-0.9);
  \draw[target] (\xC,1.5) -- (\xE,1.5)  (\xC,-1.5) -- (\xE,-1.5);
  \node[font=\footnotesize,anchor=south west,fill=white,inner sep=1pt] at (\xA+0.05,0.93) {$\pm\log\log\log t$};
  \node[font=\footnotesize,anchor=north] at ({(\xC+\xD)/2},1.45) {$\pm\log\log\log(t+1)$};
  \draw[bpath] (0,0.3) -- (\xA,0.3) -- (\xH,{-f(\xH)})
      -- plot[domain=\xH:\xB,samples=10] (\x,{-f(\x)})
      -- (\xC,1.026) -- (\xD,-1.934) -- (\xK,{f(\xK)})
      -- plot[domain=\xK:\xE,samples=10] (\x,{f(\x)});
  \fill (\xB,{-f(\xB)}) circle (1.6pt);
  \node[font=\scriptsize,anchor=north,align=center] at (\xB,{-f(\xB)-0.12})
       {$\dfrac{b}{a}\le\dfrac1{\log\log t}$};
  \fill (\xK,{f(\xK)}) circle (1.6pt);
  \node[font=\scriptsize,anchor=south east] at (\xK-0.1,{f(\xK)+0.08}) {attains $a\log n$};
\end{scope}
\end{tikzpicture}}
\usepackage{url}
\usepackage{amsfonts,amstext,amsbsy,amssymb, amsmath}
\usepackage{amsopn,amsthm,amscd,amsxtra}
\usepackage{latexsym,mathrsfs} 
\usepackage{mathabx} 
\savesymbol{widering} 
\restoresymbol{widering}{MABwidering}

\usepackage{hyperref}
\RequirePackage[dvipsnames]{xcolor} 
\definecolor{halfgray}{gray}{0.55} 
\definecolor{webgreen}{rgb}{0,0.5,0}
\definecolor{webbrown}{rgb}{.6,0,0} \hypersetup{%
  colorlinks=true, linktocpage=true, pdfstartpage=3,
  pdfstartview=FitV,%
  breaklinks=true, pdfpagemode=UseNone, pageanchor=true,
  pdfpagemode=UseOutlines,%
  plainpages=false, bookmarksnumbered, bookmarksopen=true,
  bookmarksopenlevel=1,%
  hypertexnames=true,
  pdfhighlight=/O,
  urlcolor=webbrown, linkcolor=RoyalBlue,
  citecolor=webgreen, 
  pdftitle={},%
  pdfauthor={},%
  pdfsubject={2000 Mathematical Subject Classification: Primary:},%
  pdfkeywords={},%
  pdfcreator={pdfLaTeX},%
  pdfproducer={LaTeX with hyperref}%
}

\newtheorem{theorem}{Theorem}[section]
\newtheorem{lemma}[theorem]{Lemma}
\newtheorem{proposition}[theorem]{Proposition}
\newtheorem{corollary}[theorem]{Corollary}

\newtheorem{maintheorem}{Theorem}

\newtheorem{maincorollary}[maintheorem]{Corollary}
\theoremstyle{definition}
\newtheorem{definition}[theorem]{Definition}
\theoremstyle{remark}
\newtheorem{remark}[theorem]{Remark}
 
\newcommand{\N}{\mathbb{N}}
\newcommand{\Z}{\mathbb{Z}}
\newcommand{\R}{\mathbb{R}}
\newcommand{\cO}{\mathcal{O}}
\newcommand{\OO}{\mathbb{O}}
\newcommand{\Ob}{\overline{\mathbb{O}}}
\newcommand{\cB}{\mathcal{B}}
\newcommand{\PP}{\mathbb{P}}
\newcommand{\EE}{\mathbb{E}}
\newcommand{\cL}{\mathcal{L}}
\newcommand{\ord}[1]{\llbracket #1 \rrbracket}
\newcommand{\sep}{\operatorname{sep}}
\newcommand{\spn}{\operatorname{span}}
\newcommand{\htop}{h_{\mathrm{top}}}
\newcommand{\hpol}{h_{\mathrm{pol}}}
\newcommand{\pair}[2]{\langle #1,#2\rangle}
\newcommand{\floor}[1]{\lfloor #1\rfloor}
\newcommand{\ceil}[1]{\lceil #1\rceil}
 
\begin{document}
 
\title{Minimal dynamics with polynomial and intermediate complexity}
 
\author{Heric Corr\^ea}
\address{IMPA. Estrada Dona Castorina, 110. Rio de Janeiro, 22460--320, Brazil.}
\email{heric.silva@impa.br}
\urladdr{heric-correa.github.io}
\thanks{The first author was supported by a CNPq grant.}
 
\author{Gabriel Lacerda}
\address{UFRJ. Av. Athos da Silveira Ramos, 149. Rio de Janeiro, 21941--909, Brazil.}
\email{lacerda@im.ufrj.br}
\urladdr{sites.google.com/view/lacerda}
\thanks{The second author was supported by a CNPq grant.}
 
 
\begin{abstract}
We study how complicated the orbits of a minimal system can be when the topological entropy vanishes, by looking at the exact order of growth of the number of distinguishable orbit segments, and not only at its exponent. For a broad class of growth rates, including all polynomial ones, we construct minimal homeomorphisms of one-dimensional continua realizing them; in particular, two such systems may have the same polynomial entropy and yet incomparable growth. We also show that passing from a minimal subshift to the maps obtained by running its suspension flow for a fixed time only rescales time in the growth law, and we characterize exactly when this rescaling has no effect. Combined with a construction of Cassaigne, this yields minimal homeomorphisms of continua whose growth lies strictly between all polynomial and all exponential rates.
\end{abstract}
 
\maketitle
 
\section{Introduction}
 
For a dynamical system with zero topological entropy, the number of distinguishable orbit segments may grow at many different subexponential rates. Polynomial entropy records only part of this information: the growth laws $n^2$ and $n^2\log n$, for example, both have polynomial exponent $2$. The generalized entropy introduced by Correa and Pujals \cite{CP} retains the growth order, up to multiplicative constants. Its values need not be comparable: the ratio of two growth functions may be arbitrarily large and arbitrarily small along different sequences of times.

We study which growth orders occur for minimal homeomorphisms of one-dimensional continua. The two requirements play different roles. Minimality forces every orbit to be dense, while connectedness excludes direct use of symbolic examples on Cantor sets. Neither condition should be confused with low complexity: minimal homeomorphisms can have positive topological entropy, as shown by Rees \cite{Rees}. Our concern here is the finer range of possibilities when topological entropy vanishes.

The construction begins with minimal subshifts of prescribed word complexity and passes to suitable time maps of their suspension flows. This produces connected phase spaces. To retain the prescribed growth order, however, one must control the change of time scale introduced by the suspension. Theorem~\ref{thm:C} gives that control, and the realization results follow by combining it with symbolic constructions and growth functions with oscillating ratios.

Write $\ord{a(n)}$ for the growth order of a nondecreasing sequence $a$. The partially ordered set of growth orders is denoted by $\OO$, and its Dedekind--MacNeille completion by $\Ob$. We use the polynomial and exponential families
\[\PP=\{\ord{n^d}:d>0\},\qquad
\EE=\{\ord{e^{hn}}:h>0\}.\]
Their supremum and infimum, respectively, belong to $\Ob\setminus\OO$. Thus the inequalities $\sup\PP<o(f)<\inf\EE$ locate the generalized entropy strictly between the polynomial and exponential scales. Definitions and the numerical projections giving polynomial and topological entropy are recalled in Section~\ref{sec:orders}.

\subsection{Main results}

Let $\cO$ be the set of nondecreasing maps $a\colon\N\to\N$, and consider the sequences with the \emph{bounded doubling property}
\begin{equation}\label{eq:B}
\cB:=\Big\{a\in\cO:\ \sup_{n\in\N}\frac{a(2n)}{a(n)}<\infty\Big\}.
\end{equation}
Typical examples are $a(n)=\ceil{n^d}$, $d\ge0$, and $a(n)=1+\floor{\log_2 n}$. For a continuous map $f$ of a compact metric space, $o(f)\in\Ob$ denotes its generalized entropy (Section~\ref{sec:orders})

\begin{maintheorem}\label{thm:A}
For every $a\in\cB$ there exists a minimal homeomorphism $f_a$ of a one-dimensional continuum $M_a$ such that
\[o(f_a)=\ord{n\,a(n)}.\]
Moreover, if $a$ is unbounded, there exists $b\in\cB$ such that $o(f_b)$ and $o(f_a)$ are incomparable and $\hpol(f_b)=\hpol(f_a)$.
\end{maintheorem}

The spaces $M_a$ are mapping tori of infinite minimal subshifts; in particular, they are one-dimensional matchbox manifolds (Remark~\ref{rem:matchbox}). The class $\ord{na(n)}$ comes from the symbolic realization theorem used below. Its linear factor is consistent with the Morse--Hedlund bound $p_\Lambda(n)\ge n+1$, but that bound alone does not imply the form $na(n)$ with $a\in\cB$. Bounded doubling ensures that the growth order is preserved under the time rescaling in Theorem~\ref{thm:C}.

\begin{maincorollary}\label{cor:B}
For every $\alpha\in[1,\infty)$ there exist two minimal homeomorphisms of one-dimensional continua with polynomial entropy $\alpha$ and mutually incomparable generalized entropies.
\end{maincorollary}

K\r{u}rka constructed Toeplitz subshifts with prescribed lower and upper polynomial entropies \cite[Proposition~4.79]{Kurka}; these minimal systems are defined on Cantor sets. Roth, Roth and Snoha realized every polynomial entropy in $[0,\infty]$ by homeomorphisms of continua \cite[Theorem~38]{RRS}. Their cone construction has a fixed vertex and therefore does not give minimal homeomorphisms on nondegenerate continua. Here we impose minimality and connectedness simultaneously and retain the full growth order.

The symbolic realization in Theorem~\ref{thm:symbolic} is a bounded-doubling case of the result attributed to Goyon in \cite[p.~149]{Ferenczi} and \cite[Theorem~1.2.35]{Nekrashevych}. We give a proof from K\r{u}rka's construction, including the quantitative bounds needed here. The passage to continua uses the suspension formula below. Incomparability at a fixed polynomial entropy is obtained from Lemma~\ref{lem:incomparable}, which preserves both bounded doubling and the polynomial exponent.

The suspension formula does not require minimality of the underlying subshift.

\begin{maintheorem}\label{thm:C}
Let $\sigma\colon\Lambda\to\Lambda$ be a nonempty two-sided subshift over a finite alphabet, and let $\Phi$ be its suspension flow with constant roof $1$. Then, for every $t\ne0$,
\[o(\Phi^t)=\ord{p_\Lambda(\ceil{|t|n})}.\]
In particular:
\begin{enumerate}
\item[(i)] $o(\Phi^t)$ does not depend on $t\ne0$ if and only if $p_\Lambda\in\cB$;
\item[(ii)] if $\sup\PP<o(\sigma)<\inf\EE$, then $\sup\PP<o(\Phi^t)<\inf\EE$ for every $t\ne0$.
\end{enumerate}
\end{maintheorem}

 The formula is the growth-order counterpart of the usual entropy identities for time maps. It also shows that equality of the generalized entropies of the time-$1$ and time-$2$ maps already implies invariance under every nonzero time rescaling. In the realization results, minimality is used to obtain a continuum and to select a minimal time map.

For growth above every polynomial, bounded doubling is unavailable. The same formula nevertheless preserves the intermediate range and gives the following applications.

\begin{maintheorem}\label{thm:D}
 The following hold:
\begin{enumerate}
\item[(a)] For every $\alpha\in(0,1)$ there exist a minimal homeomorphism $f_\alpha$ of a one-dimensional continuum and constants $0<c_\alpha<C_\alpha$ such that
\[
\ord{e^{c_\alpha n^\alpha}}\le o(f_\alpha)\le\ord{e^{C_\alpha n^\alpha}}.
\]
In particular, $\sup\PP<o(f_\alpha)<o(f_{\alpha'})<\inf\EE$ whenever $0<\alpha<\alpha'<1$.
\item[(b)] There exist minimal homeomorphisms $f$ and $g$ of one-dimensional continua such that $\sup\PP<o(f)<\inf\EE$, $\sup\PP<o(g)<\inf\EE$, and $o(f)$ and $o(g)$ are incomparable.
\end{enumerate}
\end{maintheorem}

For Theorem~\ref{thm:D} we use Cassaigne's uniformly recurrent words with $\log p(n)\sim\varphi(n)$ \cite[Theorem~3]{Cassaigne}. This construction includes both $\varphi(n)=n^\alpha$ and the oscillating function used in part~(b). The suspension formula transfers the resulting estimates to continua; comparison of the oscillating example with $\varphi(n)=n^{1/2}$ gives the incomparable pair. Part~(a) gives stretched-exponential bounds, rather than the exact class $\ord{e^{n^\alpha}}$.

Suspensions of Cantor systems are a standard source of connected minimal examples. Boro\'nski, Clark and Oprocha \cite{BCO} develop pseudo-suspensions to obtain continua with additional topological and dynamical properties. The estimates considered here concern generalized entropy in the zero-topological-entropy regime. They also provide examples for the realization problem formulated by Correa and Pujals \cite[p.~1592]{CP}, now within minimal homeomorphisms of one-dimensional continua.

\subsection*{Organization of the paper.} Sections~\ref{sec:prelim}--\ref{sec:orders} contain preliminaries, orders of growth and the class $\cB$; Section~\ref{sec:symbolic} the symbolic constructions; Section~\ref{sec:continua} the analysis of suspensions; Section~\ref{sec:proofs} the proofs of the main results; and Section~\ref{sec:questions} some questions.

 
\section{Preliminaries}\label{sec:prelim}

A continuum is a nonempty compact connected metric space. Throughout, dimension means covering dimension. We write $\N=\{1,2,\ldots\}$ for lengths and growth sequences; one-sided words are indexed by $\Z_{\ge0}$. Unless indicated otherwise, $\log=\log_2$, while $\ln$ denotes the natural logarithm.
 
\subsection{Suspensions}
Let $f$ be a homeomorphism of a compact metric space $X$. The \emph{mapping torus} of $f$ is \[M_X:=(X\times[0,1])/\!\sim,\] where $(x,1)\sim(f(x),0)$; we write $[x,s]$ for the class of $(x,s)$. The \emph{suspension flow} of $f$ is $\phi^t([x,s]):=[x,s+t]$, with the identifications understood. A flow is \emph{minimal} if every orbit is dense. We shall use the following result to pass from minimal flows to minimal homeomorphisms.
 
\begin{proposition}[{\cite[Proposition~1.5]{IO}}]\label{prop:timet}
If $\phi$ is a minimal flow on a compact metric space, then $\phi^t$ is minimal for every $t\in\R$ outside a countable set.
\end{proposition}
 
\subsection{Subshifts}

Let $A$ be a finite alphabet and let $\Sigma=A^\Z$ (or $A^{\Z_{\ge0}}$), with shift \[\sigma(\xi)=(\xi_{n+1})_n.\] A subshift is the restriction of $\sigma$ to a nonempty closed set $\Lambda$ satisfying $\sigma(\Lambda)=\Lambda$. It is minimal if every forward orbit is dense. For $r\le s$, write $\xi_{[r,s]}=(\xi_r,\ldots,\xi_s)$. The language of length $n$ is the set $\cL_n(\Lambda)$ of all words $\xi_{[i,i+n-1]}$ occurring in points of $\Lambda$, and $\cL(\Lambda)=\bigcup_{n\ge1}\cL_n(\Lambda)$. The complexity function $p_\Lambda(n)=\#\cL_n(\Lambda)$ is nondecreasing.

A sequence $\xi$ is \emph{Toeplitz} if every coordinate is periodic: for every $n$ there is $\tau(n)\ge1$ with $\xi_{n+k\tau(n)}=\xi_n$ for all $k$. A \emph{Toeplitz shift} is the orbit closure of a nonperiodic Toeplitz sequence; since Toeplitz sequences are uniformly recurrent, Toeplitz shifts are minimal \cite[p.~191]{Kurka}.
 
\section{Orders of growth and generalized entropy}\label{sec:orders}
 
\subsection{Basic definitions}
For $a,b\in\cO$ write $a\preccurlyeq b$ if $a(n)\le Cb(n)$ for some $C>0$ and all $n$, and $a\asymp b$ if $a\preccurlyeq b\preccurlyeq a$. The class of $a$ is $\ord{a(n)}$, and $\OO=\cO/\!\asymp$ is ordered by $\ord{a(n)}\le\ord{b(n)}$ iff $a\preccurlyeq b$. As in \cite{CP}, we work in the Dedekind--MacNeille completion $\Ob\supset\OO$, a complete lattice. Two elements of $\Ob$ are \emph{incomparable} if neither is below the other; for $\ord{a_1(n)},\ord{a_2(n)}\in\OO$ this means
\[
\limsup_{n\to\infty}\frac{a_1(n)}{a_2(n)}=\infty\quad\text{and}\quad\liminf_{n\to\infty}\frac{a_1(n)}{a_2(n)}=0.
\]
 
Let $f$ be a continuous map of a compact metric space $(X,d)$. For $\varepsilon>0$ and $n\in\N$, let $B^n_f(x,\varepsilon):=\{y: \max_{0\le k<n}d(f^kx,f^ky)<\varepsilon\}$. A set $E$ is \emph{$(n,\varepsilon)$-separated} if $E\cap B^n_f(x,\varepsilon)=\{x\}$ for every $x\in E$, and \emph{$(n,\varepsilon)$-spanning} if $X=\bigcup_{x\in E}B^n_f(x,\varepsilon)$. Let $\sep_{f,\varepsilon}(n)$ be the maximal cardinality of an $(n,\varepsilon)$-separated set and $\spn_{f,\varepsilon}(n)$ the minimal cardinality of an $(n,\varepsilon)$-spanning set. The \emph{generalized entropy} of $f$ is
\begin{equation}\label{eq:o}
o(f):=\sup_{\varepsilon>0}\ord{\sep_{f,\varepsilon}(n)}=\sup_{\varepsilon>0}\ord{\spn_{f,\varepsilon}(n)}\in\Ob.
\end{equation}
Correa and Pujals \cite{CP} define the projections $\pi_\EE,\pi_\PP\colon\Ob\to[0,\infty]$,
\[
\pi_\EE(o)=\inf\{\lambda>0:\ord{e^{\lambda n}}\ge o\},\qquad \pi_\PP(o)=\inf\{d>0:\ord{n^d}\ge o\},
\]
and prove that $\htop(f)=\pi_\EE(o(f))$ \cite[Theorem~2]{CP} and, by the same argument, $\hpol(f)=\pi_\PP(o(f))$ \cite[Section~1]{CP}, where $\hpol$ is the polynomial entropy introduced by Marco \cite{Marco}; see also \cite{Marco09}.
 
\begin{remark}\label{rem:facts}
We record some facts used repeatedly. (i) Sequences in $\cO$ take values in $\N=\{1,2,\dots\}$; a real nondecreasing sequence $c\ge1$ is identified with $\ceil{c}\in\cO$, which does not affect orders since $c\le\ceil{c}\le2c$. (ii) If $a(n)/b(n)\to0$, then $\ord{a(n)}<\ord{b(n)}$. (iii) $o(f)$ is invariant under topological conjugacy and, for a homeomorphism, $o(f^{-1})=o(f)$ \cite[Theorem~1 and Proposition~B.2]{CP}. (iv) If $c\in\cO$, then $\pi_\PP(\ord{c(n)})=\limsup_n \log c(n)/\log n$: if $d$ exceeds the $\limsup$ then $c(n)\le n^d$ for large $n$, and if $d$ is smaller then $c(n)/n^d$ is unbounded.
\end{remark}
 
\subsection{Bounded doubling and affine invariance}
 
The role of $\cB$ is to make orders of growth insensitive to the time rescalings $n\mapsto\ceil{An}+B$; for instance $e^{\sqrt{2n}}/e^{\sqrt n}\to\infty$, so $\ord{e^{\sqrt{2n}}}>\ord{e^{\sqrt n}}$, whereas polynomial orders are unaffected.
 
Given $a\in\cO$, $A>0$ and $B\in\Z$, set $a_{A,B}(n):=a(\max\{1,\ceil{An}+B\})$.
 
\begin{definition}
A sequence $a\in\cO$ is \emph{invariant under affine transformations} if $\ord{a(n)}=\ord{a_{A,B}(n)}$ for all $A>0$ and $B\in\Z$.
\end{definition}
 
\begin{lemma}\label{lem:affine}
If $a\in\cB$, then $n\,a(n)\in\cB$. Every sequence in $\cB$ is invariant under affine transformations.
\end{lemma}
 
The second assertion is a standard property of O-\hspace{0pt}regularly varying sequences (see Remark~\ref{rem:B} and \cite[Chapter~2]{BGT}); we include the short proof for completeness.
 
\begin{proof}
If $a(2n)\le Ca(n)$, then $b(n):=na(n)$ satisfies $b(2n)\le 2Cb(n)$. Now let $c\in\cB$ with $c(2n)\le Cc(n)$, $C\ge1$, fix $A>0$, $B\in\Z$, and let $\ell_n:=\max\{1,\ceil{An}+B\}$. Since $An+B\le\ceil{An}+B<An+B+1$, there are $m,N\in\N$ with $2^{-m}n\le\ell_n\le 2^mn$ for $n\ge N$. By monotonicity and $m$ applications of doubling, $c(\ell_n)\le c(2^mn)\le C^mc(n)$ and $c(n)\le c(2^m\ell_n)\le C^mc(\ell_n)$ for $n\ge N$. Since $c\ge1$, the finitely many $n<N$ only change the constants, so $\ord{c_{A,B}(n)}=\ord{c(n)}$.
\end{proof}
 
Conversely, $A=2$, $B=0$ shows that affine invariance implies bounded doubling, so within $\cO$ the two notions coincide.
 
\begin{remark}\label{rem:B}
(i) A nondecreasing sequence $a$ belongs to $\cB$ if and only if $\limsup_n a(\ceil{\lambda n})/a(n)<\infty$ for some, equivalently every, $\lambda>1$. Thus $\cB$ is exactly the class of nondecreasing \emph{O-\hspace{0pt}regularly varying} sequences in the sense of Karamata; see \cite[Chapter~2]{BGT} for O-\hspace{0pt}regular variation and \cite[Section~3]{JGSS} for sequences.
 
(ii) Iterating $a(2n)\le Ca(n)$ gives $a(n)\le Ca(1)\,n^{\log_2C}$. Hence every element of $\cB$ is bounded by a polynomial, $\ord{a(n)}<\sup\PP$ for every $a\in\cB$, and no order above $\sup\PP$ has the bounded doubling property.
\end{remark}
 
We now construct, for $a\in\cB$, a sequence $b\in\cB$ that is infinitely often much larger and infinitely often much smaller than $a$, with the same polynomial exponent. Since $b$ must stay in $\cB$, its oscillations have to be built slowly.
 
\begin{lemma}\label{lem:incomparable}
Given $a\in\cB$ unbounded, there exists $b\in\cB$ such that
\begin{equation}\label{eq:incomp}
\limsup_{n\to\infty}\frac{b(n)}{a(n)}=\infty\qquad\text{and}\qquad\liminf_{n\to\infty}\frac{b(n)}{a(n)}=0.
\end{equation}
Moreover,
\begin{equation}\label{eq:samepol}
\limsup_{n\to\infty}\frac{\log a(n)}{\log n}=\limsup_{n\to\infty}\frac{\log b(n)}{\log n}.
\end{equation}
\end{lemma}
 
\begin{figure}[ht]
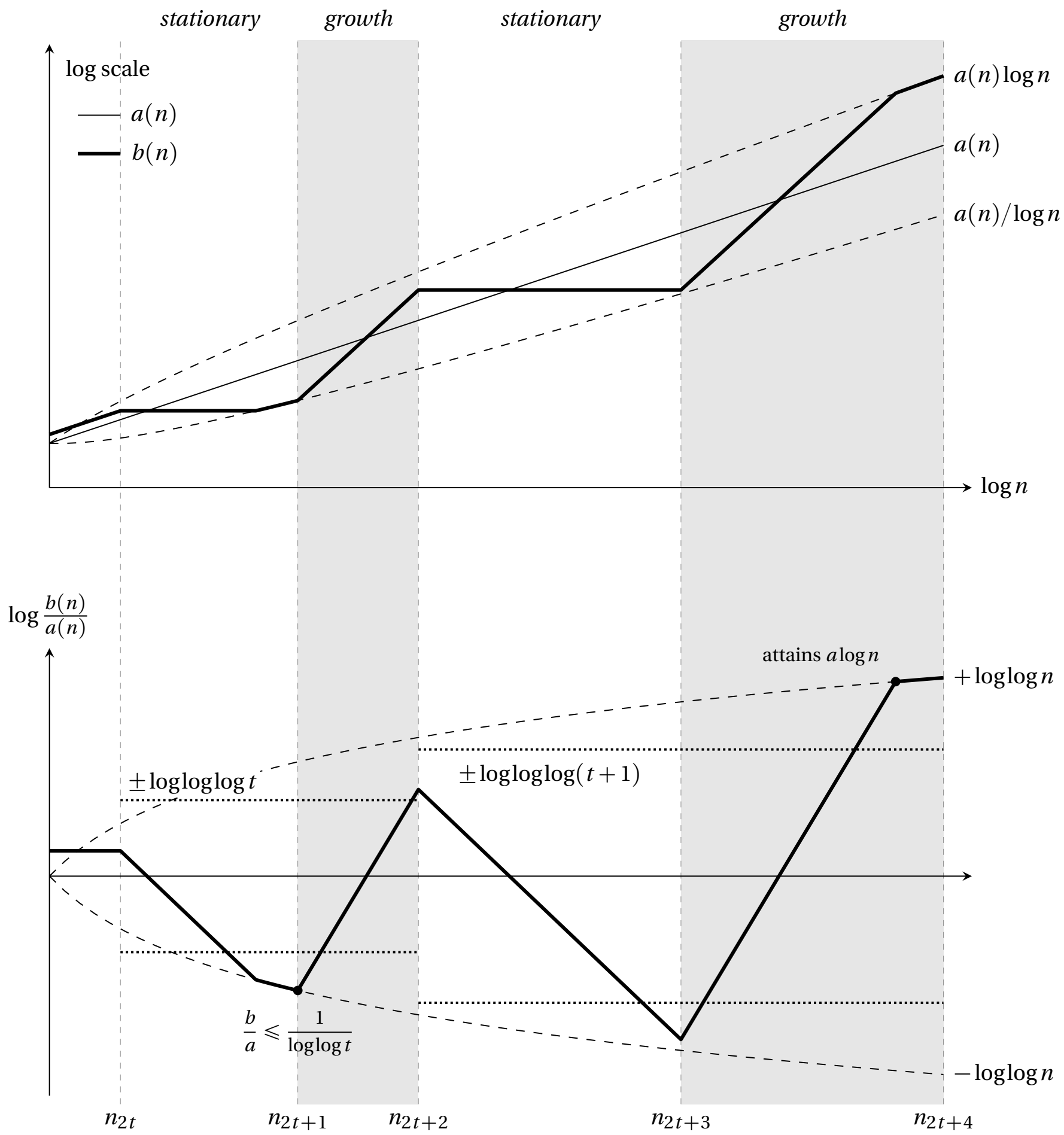

\centering
\resizebox{0.94\textwidth}{!}{\graficosdaprova}
\caption{Schematic of the construction in logarithmic scale. The lower panel shows $\log(b/a)$ and its logarithmic envelopes.}\label{fig:construction}
\end{figure}

\noindent\emph{Construction strategy.}
We alternate phases in which $b$ is kept as small as the lower envelope allows with phases in which it grows at the maximal permitted doubling rate. Throughout the construction,
\[
\frac{a(n)}{\log n}\le b(n)\le a(n)\log n,
\qquad b(2n)\le Kb(n).
\]
The envelope preserves the polynomial exponent, while the alternating phases make $b/a$ arbitrarily small and arbitrarily large. Figure~\ref{fig:construction} illustrates the construction; the recursive definition below also handles flat portions of $a$. 

\begin{proof}
Throughout, $\log=\log_2$. Let $C\ge1$ satisfy $a(2n)\le Ca(n)$, and fix an integer $K>C$. We construct a nondecreasing $b\colon\N\to[1,\infty)$ and integers $n_8<n_9<\cdots$ such that
\begin{enumerate}
\item[(i)] $b(n)\le Kb(\ceil{n/2})$ and (ii) $a(n)/\log n\le b(n)\le a(n)\log n$ for every $n\ge2$;
\item[(iii)] for every $t\ge4$, $b(n_{2t+1})\le a(n_{2t+1})/\log\log t$, and $b(n'_t)\ge a(n'_t)\log\log t$ for some $n'_t\in(n_{2t+1},n_{2t+2}]$.
\end{enumerate}
Put $b:=a$ on $[1,n_8]$, with $n_8\ge16$; then (i) and (ii) hold there. The key observation is that the lower envelope doubles by at most $C<K$: if $n\ge3$, $m:=\ceil{n/2}$ and (ii) holds at $m$, then $a(n)\le a(2m)\le Ca(m)$ and $\log n\ge\log m$, so
\begin{equation}\label{eq:env}
a(n)/\log n\le Cb(m)\le Kb(\ceil{n/2}).
\end{equation}
Assume $b$ is defined on $[1,n_{2t}]$ with (i) and (ii). Since $n_{2t}\ge2t$, we have $\log n\ge\log\log t$ for $n>n_{2t}$.
 
\emph{Stationary phase.} Choose $n_{2t+1}>n_{2t}$ with $a(n_{2t+1})\ge b(n_{2t})\log\log t$ and $\log n_{2t+1}\ge\log\log t$ (possible since $a$ is unbounded), and put \[b(n):=\max\{b(n-1),a(n)/\log n\} \quad \text{for} \ n \in (n_{2t},n_{2t+1}]. \] Then (i) follows from (i) at $n-1$ and \eqref{eq:env}, and (ii) from the monotonicity of $a$ and $\log$. Since $a(n)/\log n\le a(n_{2t+1})/\log\log t$ on this interval, we get $b(n_{2t+1})\le a(n_{2t+1})/\log\log t$.
 
\emph{Growth phase.} Choose $s\ge1$ with \[(K/C)^s\,b(n_{2t+1})/a(n_{2t+1})\ge\log\log t,\] put $n_{2t+2}:=2^sn_{2t+1}$ and $b(n):=\min\{Kb(\ceil{n/2}),a(n)\log n\}$ on $(n_{2t+1},n_{2t+2}]$. Property (i) and the upper bound in (ii) hold by definition, the lower bound by \eqref{eq:env}, and monotonicity follows from (i) and (ii) at $n-1$. If $b(n')=a(n')\log n'$ for some $n'$ in this interval, then $b(n')/a(n')=\log n'\ge\log n_{2t+1}\ge\log\log t$. Otherwise $b(n_{2t+2})=K^sb(n_{2t+1})$, while $a(n_{2t+2})\le C^sa(n_{2t+1})$, so $b(n_{2t+2})/a(n_{2t+2})\ge\log\log t$. This proves (iii).
 
Finally, $\tilde b:=\floor{b}\in\cO$ satisfies $b/2\le\tilde b\le b$. By (i), $\tilde b(2n)\le Kb(n)\le2K\tilde b(n)$, so $\tilde b\in\cB$; (iii) gives \eqref{eq:incomp}; and (ii) gives $|\log\tilde b(n)-\log a(n)|\le1+\log\log n$, which after division by $\log n$ yields \eqref{eq:samepol}.
 
\end{proof}
 
\subsection{Generalized entropy of a subshift}\label{sec:subshift}
 
The following argument is standard. Equip $A^\Z$ with the metric $d(\xi,\eta)=2^{-N(\xi,\eta)}$, where $N(\xi,\eta)=\min\{|j|:\xi_j\ne\eta_j\}$. Choosing one point for each word of $\cL_n(\Lambda)$, realized on the coordinates $[0,n-1]$, gives an $(n,1/2)$-separated set, since two distinct such words differ at a coordinate that one of the first $n$ iterates moves to $0$. Conversely, for $r\ge0$, choosing one point for each word of length $n+2r$ realized on $[-r,n+r-1]$ gives an $(n,2^{-r})$-spanning set. Hence
\[
\sep_{\sigma,1/2}(n)\ge p_\Lambda(n),\qquad \spn_{\sigma,2^{-r}}(n)\le p_\Lambda(n+2r)\le|A|^{2r}p_\Lambda(n),
\]
and therefore
\begin{equation}\label{eq:subshift}
o(\sigma)=\ord{p_\Lambda(n)}.
\end{equation}
 
\section{Symbolic realization}\label{sec:symbolic}
 
\subsection{K\r{u}rka's construction}
We use the following quantitative form of K\r{u}rka's construction \cite[Proposition~4.79, pp.~199--200]{Kurka}. The monotonicity assumption $q_k\le q_{k+1}$ ensures that every letter at a given level occurs in the limiting word. We include the construction and the counting argument.
 
\begin{proposition}\label{prop:kurka}
Let $(q_k)_{k\ge0}$ be integers such that, for each $k\ge0$,
\begin{equation}\label{eq:qk}
3\le q_k{}\le q_{k+1}\le(q_k-1)^2.
\end{equation}
Then there exists an infinite minimal one-sided Toeplitz shift $\sigma\colon\Lambda\to\Lambda$ such that
\begin{equation}\label{eq:kurka}
2\cdot5^{k-1}q_k\le p_\Lambda(5^k)\le2\cdot5^kq_k\qquad\text{for each } k\ge1.
\end{equation}
\end{proposition}
 
\begin{proof}

Let $A_k=\{0,1,\ldots,q_k-1\}$. For each $k\ge0$, list the first $q_{k+1}$ pairs $(b,c)\in\{1,\ldots,q_k-1\}^2$ in lexicographic order, and assign them in order to $a\in A_{k+1}$. Set
\[
\vartheta_k(a)=0\,0\,b\,0\,c,
\qquad \Theta_k=\vartheta_0\circ\cdots\circ\vartheta_{k-1},
\]
with $\Theta_0$ the identity. The inequality $q_{k+1}\le(q_k-1)^2$ makes this possible, and each $\Theta_k$ is injective on letters.

Since every $\vartheta_k(0)$ begins with $0$, the words $\Theta_k(0)$ are nested prefixes. Let $u$ be their one-sided limit and let $\Lambda$ be its shift orbit closure. More generally, put
\[
u^{(k)}=\lim_{m\to\infty}
\vartheta_k\circ\cdots\circ\vartheta_{m-1}(0).
\]
Then $u=\Theta_k(u^{(k)})$ and $u^{(k)}=\vartheta_k(u^{(k+1)})$.

Every letter of $A_k$ occurs in $u^{(k)}$. Indeed, for $1\le a\le q_k-1$, the last letter of $\vartheta_k(a-1)$ is $a$. Since $q_{k+1}\ge q_k$, we can repeat this observation at the next levels, decreasing the letter by one each time. Thus $a$ occurs in
$\vartheta_k\circ\cdots\circ\vartheta_{k+a-1}(0)$, a prefix of $u^{(k)}$; the letter $0$ already occurs at its beginning.

The word $u$ is Toeplitz. Given a coordinate $i\ge0$, choose $k$ with $i<5^k$. The first level-$k$ block in every level-$(k+1)$ block is $\Theta_k(0)$, so
\[
u_{i+\ell5^{k+1}}=u_i\qquad(\ell\ge0).
\]
Moreover, every finite factor of $u$ lies in a prefix $\Theta_k(0)$, and this prefix occurs at every multiple of $5^{k+1}$. Hence $u$ is uniformly recurrent and $\Lambda$ is minimal. It is infinite by the lower bound proved below.

For each $k$, the level-$k$ letters in $u$ form the concatenation $u^{(k)}=\vartheta_k(u^{(k+1)})$. Their zero/nonzero pattern is therefore the periodic pattern $00{*}0{*}$, where ${*}$ denotes a nonzero letter. In particular, every pair of consecutive level-$k$ letters contains a zero.

\emph{Upper bound.} A word $u\in\cL_{5^k}(\Lambda)$ is a factor of $\Theta_k(ab)$ for two consecutive letters $a,b\in A_k$ with $a=0$ or $b=0$. There are at most $2q_k$ such pairs $(a,b)$, and $5^k$ possible positions of $u$ inside $\Theta_k(ab)$. Hence $p_\Lambda(5^k)\le2\cdot5^kq_k$.
 
\emph{Lower bound.} We first explain how a word of length $5^k$ determines its alignment with the level-$k$ blocks. Assume its alignment is already known modulo $5^j$, starting with $j=0$. The word contains at least four complete level-$j$ blocks. By injectivity of $\Theta_j$, their letters, and hence their zero/nonzero indicators, can be recovered. The five possible four-letter windows of the indicator pattern $00101$ are
\[
0010,\quad0101,\quad1010,\quad0100,\quad1001.
\]
They are distinct. Thus these four blocks determine the position of their first block modulo $5$, and hence the alignment of the original word modulo $5^{j+1}$. Induction gives its alignment modulo $5^k$.

For $a\in A_k$, choose a letter $a'$ that follows $a$ in $u^{(k)}$. For $0\le s<2\cdot5^{k-1}$, let $w_{a,s}$ be the length-$5^k$ word beginning at position $s$ in $\Theta_k(a)\Theta_k(a')$. Write $\vartheta_{k-1}(a)=00b_a0c_a$. The word $w_{a,s}$ contains the complete blocks $\Theta_{k-1}(b_a)$ and $\Theta_{k-1}(c_a)$. Its alignment determines $s$, since $0\le s<5^k$. Once this alignment is known, injectivity of $\Theta_{k-1}$ recovers $b_a,c_a$, and injectivity of $\vartheta_{k-1}$ recovers $a$. Consequently, the $2\cdot5^{k-1}q_k$ words $w_{a,s}$ are distinct, proving the lower bound.
\end{proof}
 
In Proposition~\ref{prop:kurka}, $q_k$ prescribes the number of words of length $5^k$ up to the factor $5^k$, and \eqref{eq:qk} allows at most a squaring from one scale to the next. Minimality is automatic, so our symbolic constructions reduce to choosing $(q_k)$.
 
We pass to two-sided shifts by the natural extension. For a one-sided subshift $\Lambda\subset A^{\Z_{\ge0}}$, let $\hat\Lambda:=\{\xi\in A^\Z:\xi_{[r,s]}\in\cL(\Lambda)\text{ for all } r\le s\}$.
 
\begin{lemma}\label{lem:twosided}
If $\sigma\colon\Lambda\to\Lambda$ is a minimal one-sided subshift, then $\sigma\colon\hat\Lambda\to\hat\Lambda$ is a minimal two-sided subshift and $p_{\hat\Lambda}=p_\Lambda$.
\end{lemma}
 
\begin{proof}
Since $\sigma(\Lambda)=\Lambda$, every word of $\cL(\Lambda)$ extends to the left inside $\cL(\Lambda)$, and a compactness argument gives $\cL_n(\hat\Lambda)=\cL_n(\Lambda)$ for all $n$. Given $u\in\cL(\Lambda)$, minimality and compactness give $N$ such that the sets $\sigma^{-j}[u]$, $0\le j\le N$, cover $\Lambda$, where $[u]$ is the cylinder of $u$; hence $u$ is a factor of every word of $\cL_{N+|u|}(\Lambda)=\cL_{N+|u|}(\hat\Lambda)$. Now let $Y\subset\hat\Lambda$ be nonempty, closed and invariant, and $y\in Y$. Every $u\in\cL(\hat\Lambda)$ is a factor of $y_{[0,N+|u|-1]}\in\cL(\hat\Lambda)$, so $\cL(Y)=\cL(\hat\Lambda)$. Since $Y$ is closed and invariant, every point all of whose words lie in $\cL(Y)$ belongs to $Y$; hence $Y=\hat\Lambda$.
\end{proof}
 
Consequently, under the hypothesis of Proposition~\ref{prop:kurka} there is an infinite minimal two-sided subshift satisfying \eqref{eq:kurka}. If $k\ge1$ and $5^k\le n<5^{k+1}$, then monotonicity and \eqref{eq:kurka} give $p_\Lambda(n)\ge p_\Lambda(5^k)\ge2\cdot5^{k-1}q_k$ and $p_\Lambda(n)\le p_\Lambda(5^{k+1})\le 10\cdot5^kq_{k+1}$, that is,
\begin{equation}\label{eq:interp}
\tfrac{2}{25}\,n\,q_k\le p_\Lambda(n)\le10\,n\,q_{k+1}.
\end{equation}
 
\subsection{Polynomial-type complexity}
The following realization is a bounded-doubling case of the result attributed to Goyon in \cite[p.~149]{Ferenczi} and \cite[Theorem~1.2.35]{Nekrashevych}. We derive it from Proposition~\ref{prop:kurka}.

\begin{theorem}\label{thm:symbolic}
For every $a\in\cB$ there exists an infinite minimal two-sided subshift $\sigma\colon\Lambda\to\Lambda$ such that $\ord{p_\Lambda(n)}=\ord{n\,a(n)}$.
\end{theorem}
 
\begin{proof}
Let $C\ge1$ satisfy $a(2n)\le Ca(n)$. Monotonicity and three applications of doubling give $a(5n)\le a(8n)\le Da(n)$, with $D:=C^3$. If $a$ is bounded, it is eventually constant (it is integer-valued and nondecreasing), and we choose $k_0\ge1$ such that $5^{k_0}$ lies in the constant tail; if $a$ is unbounded, we choose $k_0\ge1$ with $a(5^{k_0})\ge D$. Let $q_k:=a(5^{k+k_0})+2\ge3$. In the bounded case $q_{k+1}=q_k\le(q_k-1)^2$. Since $a$ is nondecreasing, $q_k\le q_{k+1}$ in all cases. In the unbounded case, writing $u=a(5^{k+k_0})\ge D$, we get $q_{k+1}\le Du+2\le u^2+2\le(u+1)^2=(q_k-1)^2$. Thus \eqref{eq:qk} holds and \eqref{eq:interp} applies.
 
Let $n\ge5$ and $5^k\le n<5^{k+1}$, $k\ge1$. Since $k_0\ge1$, $q_k\ge a(5^{k+1})\ge a(n)$, so \eqref{eq:interp} gives $p_\Lambda(n)\ge\frac{2}{25}na(n)$. On the other hand, \[q_{k+1}=a(5^{k_0+1}5^k)+2\le D^{k_0+1}a(5^k)+2\le(D^{k_0+1}+2)a(n),\] so $p_\Lambda(n)\le10(D^{k_0+1}+2)\,na(n)$. Enlarging the constants to cover $1\le n<5$, we obtain \[p_\Lambda\asymp na(n).\]
\end{proof}
 
\subsection{Intermediate complexity}
 
For orders above $\sup\PP$ we use the following construction of Cassaigne. For a word $u\in\{0,1\}^{\Z_{\ge0}}$, $p_u(n)$ denotes the number of distinct factors of length $n$ of $u$, and $\log=\log_2$.
 
\begin{theorem}[{\cite[Theorem~3]{Cassaigne}}]\label{thm:cassaigne}
Let $\varphi\colon(0,\infty)\to(0,\infty)$ be differentiable, with $\varphi'$ decreasing, and assume that $\log t/\varphi(t)\to0$ and $t^{\beta}\varphi'(t)\to0$ as $t\to\infty$, for some $\beta>0$. Then there exists a uniformly recurrent word $u\in\{0,1\}^{\Z_{\ge0}}$ such that $\log p_u(n)\sim\varphi(n)$ as $n\to\infty$.
\end{theorem}
 
As observed in \cite{Cassaigne}, the hypotheses hold for $\varphi(t)=t^\alpha$ with $0<\alpha<1$, and for
\begin{equation}\label{eq:exotic}
\varphi_1(t)=(t+10)^{\frac12+\frac14\cos(\ln\ln(t+10))},
\end{equation}
which oscillates slowly between $t^{1/4}$ and $t^{3/4}$. Since $u$ is uniformly recurrent, its orbit closure is a minimal one-sided subshift with complexity $p_u$, and Lemma~\ref{lem:twosided} gives an infinite minimal two-sided subshift $\Lambda_\varphi$ with $p_{\Lambda_\varphi}=p_u$.
 
\section{From subshifts to continua}\label{sec:continua}
 
Throughout this section, $\sigma\colon\Lambda\to\Lambda$ is a nonempty two-sided subshift over a finite alphabet, $\Phi$ is its unit-roof suspension flow on $M_\Lambda$, and $D$ is a compatible metric on $M_\Lambda$. Minimality and infinitude will be imposed explicitly when needed. The entropy formula is independent of the choice of $D$. For $\xi\in\Lambda$ and $v\in\R$, set
\begin{equation}\label{eq:param}
\pair{\xi}{v}:=\big[\sigma^{\floor{v}}\xi,\,v-\floor{v}\big]\in M_\Lambda.
\end{equation}
On each $\Lambda\times[m,m+1]$ this coincides with $(\xi,v)\mapsto[\sigma^m\xi,v-m]$, so $(\xi,v)\mapsto\pair{\xi}{v}$ is continuous. Moreover, for $m\in\Z$, $v,t\in\R$ and $s\in[0,1)$,
\begin{equation}\label{eq:param2}
\pair{\xi}{v+m}=\pair{\sigma^m\xi}{v},\qquad\Phi^t\pair{\xi}{s}=\pair{\xi}{s+t},\qquad\pair{\xi}{s}=[\xi,s].
\end{equation}
 
For $t>0$, a point $\pair{\xi}{v}$ can be viewed as a sequence with a reading position moving at unit speed. The first $n$ iterates of $\Phi^t$ inspect a block of length about $tn$. The estimates below make this observation precise; negative times are handled by invariance under inversion.
 
\begin{proposition}\label{prop:continuum}
If $\Lambda$ is infinite and minimal, then $M_\Lambda$ is a one-dimensional continuum and $\Phi$ is a minimal flow.
\end{proposition}

\begin{proof}
    Every nonempty open set $U\subset M_\Lambda$ contains a subset of the form \[\{[y,r]: y\in V_0,\ r\in I\},\] where $V_0\subset\Lambda$ is open and $I\subset(0,1)$ is an open interval. Given $[x,s]$, the minimality of $\sigma$ yields an integer $n\ge0$ such that $\sigma^nx\in V_0$. Consequently, $\Phi^{n-s+r}[x,s]\in U$ for $r\in I$; hence $\Phi$ is minimal. $M_\Lambda$ is compact as a quotient of $\Lambda\times[0,1]$, and connected because it is the closure of a single (connected) orbit. 
    
For the dimension, $\Lambda$ is zero-dimensional as a closed subspace of the product of finite discrete spaces $A^\Z$ \cite[Theorem~1.3.6 and Proposition~1.2.2]{Engelking}. The quotient map is injective on each of $\Lambda\times[0,1/2]$ and $\Lambda\times[1/2,1]$. Their images are closed, cover $M_\Lambda$, and have dimension at most $1$ by the product theorem \cite[Theorem~3.2.13]{Engelking}. The closed sum theorem \cite[Theorem~3.1.18]{Engelking} gives $\dim M_\Lambda\le1$. A short closed orbit segment is an embedded interval, so $\dim M_\Lambda\ge1$.
\end{proof}

\begin{remark}\label{rem:matchbox}
If $\Lambda$ is infinite and minimal, then $M_\Lambda$ is a one-dimensional matchbox manifold in the sense of \cite[Definition~2.1]{CHL}, with Cantor transversal $\Lambda$. Indeed, \[M_\Lambda=(\Lambda\times\mathbb R)/\mathbb Z\] for the free, properly discontinuous action $(\xi,v)\mapsto(\sigma\xi,v-1)$, so the quotient map $P(\xi,v)=\langle\xi,v\rangle$ is open, and it is injective on $\Lambda\times(c-\tfrac12,c+\tfrac12)$, since two points of that interval differ by less than $1$. Hence $\Psi_c(\xi,s)=\langle\xi,c+s/4\rangle$ is a homeomorphism of $\Lambda\times(-2,2)$ onto an open set, sending each $\{\xi\}\times(-2,2)$ onto an arc of a $\Phi$-orbit. Its restriction to $\Lambda\times[-1,1]$ is a homeomorphism onto the closure of $U_c:=\Psi_c(\Lambda\times(-1,1))$, and every $\langle\xi_0,c\rangle$ equals $\Psi_c(\xi_0,0)$. These are the required charts and their extensions.
\end{remark}
 
\begin{proposition}\label{prop:flow}
For every $t\ne0$, \[o(\Phi^t)=\ord{p_\Lambda(\ceil{|t|n})}.\] In particular, if $p_\Lambda$ is invariant under affine transformations, then $o(\Phi^t)=\ord{p_\Lambda(n)}$ for every $t\ne0$.
\end{proposition}
 
\begin{lemma}\label{lem:metric} The following hold:
\begin{enumerate}
\item[(a)] For every $\varepsilon>0$ there exist $L\in\N$ and $\delta\in(0,1]$ such that $D(\pair{\xi}{v},\pair{\eta}{w})<\varepsilon$ whenever $|v-w|<\delta$ and $\xi_j=\eta_j$ for $|j-\floor{v}|\le L$.
\item[(b)] For every $L\in\N$ there exists $\varepsilon_L>0$ such that $D([\xi,s],[\eta,s])\ge\varepsilon_L$ whenever $s\in[0,1]$ and $\xi_{[-L,L]}\ne\eta_{[-L,L]}$.
\end{enumerate}
\end{lemma}
 
\begin{proof}
(a) The map $P(\zeta,r):=\pair{\zeta}{r}$ is uniformly continuous on $\Lambda\times[-1,2]$: there is $\rho\in(0,1]$ such that $d(\zeta,\zeta')<\rho$ and $|r-r'|<\rho$ imply $D(P(\zeta,r),P(\zeta',r'))<\varepsilon$. Choose $L$ with $2^{-(L+1)}<\rho$ and $\delta:=\rho$. 

Given $\xi,\eta,v,w$ as in the statement, let $m=\floor{v}$, $\zeta=\sigma^m\xi$, $\zeta'=\sigma^m\eta$, $r=v-m\in[0,1)$ and $r'=w-m\in[-1,2]$. Then $\zeta,\zeta'$ agree on $[-L,L]$, so $d(\zeta,\zeta')\le2^{-(L+1)}<\rho$, and $|r-r'|<\rho$. By \eqref{eq:param2}, $\pair{\xi}{v}=P(\zeta,r)$ and $\pair{\eta}{w}=P(\zeta',r')$.
 
(b) The set $K_L$ of triples $(\xi,\eta,s)\in\Lambda\times\Lambda\times[0,1]$ with $\xi_{[-L,L]}\ne\eta_{[-L,L]}$ is compact, and $F(\xi,\eta,s):=D([\xi,s],[\eta,s])$ is continuous on it. It is also positive: $[\xi,s]=[\eta,s]$ forces $\xi=\eta$, because the classes are singletons for $s\in(0,1)$, the class of $(\xi,0)$ is $\{(\xi,0),(\sigma^{-1}\xi,1)\}$, and $[\xi,1]=[\eta,1]$ gives $\sigma\xi=\sigma\eta$. If $K_L\ne\emptyset$, take $\varepsilon_L=\min_{K_L}F$; if $K_L=\emptyset$, any $\varepsilon_L>0$ works.
\end{proof}
 
\begin{lemma}\label{lem:span}
For any $t>0$ and $\varepsilon>0$ there exist integers $A_\varepsilon\ge1$ and $B_\varepsilon\ge0$ such that \[\spn_{\Phi^t,\varepsilon}(n)\le A_\varepsilon\,p_\Lambda(\ceil{tn}+B_\varepsilon) \quad \text{for all} \ n.\]
\end{lemma}
 
\begin{proof}
Let $L,\delta$ be given by Lemma~\ref{lem:metric}(a) and fix $M\in\N$ with $1/M<\delta$. Given $n$, put $N:=\ceil{tn}+2L+1$. For each $w\in\cL_N(\Lambda)$ choose $\xi^w\in\Lambda$ with \[\xi^w_{[-L,\ceil{tn}+L]}=w,\] and let $E:=\{\pair{\xi^w}{j/M}: w\in\cL_N(\Lambda),\ 0\le j<M\}$, so $\#E\le Mp_\Lambda(N)$. 

We claim that $E$ is $(n,\varepsilon)$-spanning for $\Phi^t$. Let $x=\pair{\eta}{u}$ with $u\in[0,1)$, choose $j$ with $j/M\le u<(j+1)/M$, let $w:=\eta_{[-L,\ceil{tn}+L]}$ and $y:=\pair{\xi^w}{j/M}$. For $0\le k<n$, \eqref{eq:param2} gives $\Phi^{kt}y=\pair{\xi^w}{v_k}$ and $\Phi^{kt}x=\pair{\eta}{v'_k}$, with $v_k=j/M+kt$, $v'_k=u+kt$ and $|v_k-v'_k|<\delta$. The integer $m_k:=\floor{v_k}$ satisfies $0\le m_k\le v_k<1+(n-1)t\le tn+1$, hence $0\le m_k\le\ceil{tn}$, and $\xi^w,\eta$ agree on $[m_k-L,m_k+L]$. By Lemma~\ref{lem:metric}(a), $x\in B^n_{\Phi^t}(y,\varepsilon)$. The lemma holds with $A_\varepsilon=M$ and $B_\varepsilon=2L+1$.
\end{proof}
 
\begin{lemma}\label{lem:sep}
Let $t>0$. There exist $\varepsilon_0>0$ and $B_0,N_0\in\N$ such that \[\sep_{\Phi^t,\varepsilon_0}(n)\ge p_\Lambda(\floor{tn}-B_0) \quad \text{for all} \  n\ge N_0. \]
\end{lemma}
 
\begin{proof}
Let $L=B_0:=\ceil{t}$, let $\varepsilon_0:=\varepsilon_L$ be given by Lemma~\ref{lem:metric}(b), and let $N_0$ be the least integer with $\floor{tN_0}-B_0\ge1$. Fix $n\ge N_0$ and put $N:=\floor{tn}-B_0$. For each $w\in\cL_N(\Lambda)$ choose $\xi^w\in\Lambda$ with $\xi^w_{[0,N-1]}=w$, and let $E:=\{\pair{\xi^w}{0}:w\in\cL_N(\Lambda)\}$; these are $p_\Lambda(N)$ distinct points. Let $w\ne w'$ differ at $i\in[0,N-1]$ and put $k:=\ceil{i/t}$. Since $i\le N-1<t(n-1)$, we have $0\le k\le n-1$; moreover $i\le kt<i+t$, so $m:=\floor{kt}$ satisfies $|m-i|\le L$. By \eqref{eq:param}, $\Phi^{kt}\pair{\xi^w}{0}=[\sigma^m\xi^w,r]$ and $\Phi^{kt}\pair{\xi^{w'}}{0}=[\sigma^m\xi^{w'},r]$ with $r=kt-m\in[0,1)$, and $\sigma^m\xi^w,\sigma^m\xi^{w'}$ differ at the coordinate $i-m\in[-L,L]$. By Lemma~\ref{lem:metric}(b), $E$ is $(n,\varepsilon_0)$-separated.
\end{proof}
 
\begin{proof}[Proof of Proposition~\ref{prop:flow}]
Write $p=p_\Lambda$. By Remark~\ref{rem:facts}(iii) we may assume $t>0$. Every word of length $n$ has at most $|A|$ extensions to the right, so $p(n+B)\le|A|^Bp(n)$ for every $B\ge0$. For $\varepsilon>0$, Lemma~\ref{lem:span} gives \[\spn_{\Phi^t,\varepsilon}(n)\le A_\varepsilon|A|^{B_\varepsilon}p(\ceil{tn}),\] hence $o(\Phi^t)\le\ord{p(\ceil{tn})}$.

Conversely, since $\ceil{tn}\le(\floor{tn}-B_0)+B_0+1$, Lemma~\ref{lem:sep} gives \[p(\ceil{tn})\le|A|^{B_0+1}\sep_{\Phi^t,\varepsilon_0}(n)\] for $n\ge N_0$, hence $\ord{p(\ceil{tn})}\le o(\Phi^t)$. The last assertion follows from affine invariance.
 
\end{proof}
 
\begin{corollary}\label{cor:homeo}
If $\Lambda$ is infinite and minimal, there is $t\in(0,1)$ such that $f:=\Phi^t$ is a minimal homeomorphism of $M_\Lambda$ with $o(f)=\ord{p_\Lambda(\ceil{tn})}$. If $p_\Lambda$ is invariant under affine transformations, then $o(f)=\ord{p_\Lambda(n)}$.
\end{corollary}
 
\begin{proof}
By Propositions~\ref{prop:continuum} and \ref{prop:timet}, there is $t\in(0,1)$ such that $\Phi^t$ is minimal; the rest is Proposition~\ref{prop:flow}.
\end{proof}
 
\begin{remark}\label{rem:roof}
Let $\Phi_r$ be the suspension flow of $\sigma$ under a continuous roof function $r\colon\Lambda\to(0,\infty)$, moving at unit speed. If $r\equiv c$ is constant, then $[\xi,s]\mapsto[\xi,s/c]$ conjugates $\Phi_r^t$ to $\Phi^{t/c}$, so Proposition~\ref{prop:flow} gives $o(\Phi_r^t)=\ord{p_\Lambda(\ceil{|t|n/c})}$; the same holds if $r$ is cohomologous to a constant through a continuous transfer function, since the suspension flows are then conjugate. 

For a general continuous roof, the preceding proof does not give the same conclusion. Lemmas~\ref{lem:span} and~\ref{lem:sep} rely on synchronized reading positions: sequences agreeing on a block traverse that block at the same times. With a variable roof, differences of Birkhoff sums of $r$ produce a time drift. Whether bounded doubling alone still implies $o(\Phi_r^t)=\ord{p_\Lambda(n)}$ is left as a question in Section~\ref{sec:questions}.
\end{remark}
 
\section{Proofs of the main results}\label{sec:proofs}
 
\begin{proof}[Proof of Theorem~\ref{thm:A}]
Let $a\in\cB$. By Lemma~\ref{lem:affine}, $c(n):=na(n)$ belongs to $\cB$ and is invariant under affine transformations. By Theorem~\ref{thm:symbolic}, there is an infinite minimal two-sided subshift $\Lambda_a$ with $\ord{p_{\Lambda_a}(n)}=\ord{c(n)}$, so $p_{\Lambda_a}$ is also invariant under affine transformations. By Proposition~\ref{prop:continuum} and Corollary~\ref{cor:homeo}, $M_a:=M_{\Lambda_a}$ is a one-dimensional continuum carrying a minimal homeomorphism $f_a$ with $o(f_a)=\ord{n\,a(n)}$. 

If $a$ is unbounded, let $b\in\cB$ be given by Lemma~\ref{lem:incomparable}, and apply the first part to $b$. Since $nb(n)/(na(n))=b(n)/a(n)$, \eqref{eq:incomp} shows that $o(f_a)$ and $o(f_b)$ are incomparable. Moreover, by \eqref{eq:samepol} and Remark~\ref{rem:facts}(iv),
\[\hpol(f_b)=1+\limsup_{n\to\infty}\frac{\log b(n)}{\log n}
=1+\limsup_{n\to\infty}\frac{\log a(n)}{\log n}=\hpol(f_a).\]
\end{proof}
 
\begin{proof}[Proof of Corollary~\ref{cor:B}]
Let $\alpha\in[1,\infty)$, and let $a(n):=\floor{\log_2n}+1$ if $\alpha=1$ and $a(n):=\ceil{n^{\alpha-1}}$ if $\alpha>1$. In both cases $a$ is unbounded, $a\in\cB$, and $\limsup_n\log a(n)/\log n=\alpha-1$; by \eqref{eq:samepol}, the same holds for the sequence $b$ of Lemma~\ref{lem:incomparable}. Applying Theorem~\ref{thm:A} to $a$ and to $b$ gives minimal homeomorphisms $f_a,f_b$ of one-dimensional continua with generalized entropies $\ord{na(n)}$ and $\ord{nb(n)}$, which are incomparable by \eqref{eq:incomp}. By \cite{CP} and Remark~\ref{rem:facts}(iv),
\[\hpol(f_a)=\pi_\PP\big(\ord{na(n)}\big)=1+\limsup_{n\to\infty}\frac{\log a(n)}{\log n}=\alpha,\]
and likewise $\hpol(f_b)=\alpha$. The restriction $\alpha\ge1$ is intrinsic to our method: by the Morse--Hedlund theorem \cite{MH}, $p_\Lambda(n)\ge n+1$ for every infinite subshift (compare \cite{ACM}, where it is shown that expansive homeomorphisms of infinite compact metric spaces have polynomial entropy at least $1$).
\end{proof}
 
\begin{proof}[Proof of Theorem~\ref{thm:C}]
The formula is Proposition~\ref{prop:flow}. (i) If $p_\Lambda\in\cB$, then $p_\Lambda$ is invariant under affine transformations (Lemma~\ref{lem:affine}), so $o(\Phi^t)=\ord{p_\Lambda(n)}$ for every $t\ne0$. Conversely, if $o(\Phi^1)=o(\Phi^2)$, then $\ord{p_\Lambda(n)}=\ord{p_\Lambda(2n)}$, that is, $p_\Lambda\in\cB$. (ii) Let $t>0$. If $\ord{n^d}\le\ord{p_\Lambda(n)}$, then $\ord{n^d}=\ord{\ceil{tn}^d}\le\ord{p_\Lambda(\ceil{tn})}$; hence $o(\Phi^t)\ge\sup\PP$. If $\ord{p_\Lambda(n)}\le\ord{e^{hn}}$ for every $h>0$, then $\ord{p_\Lambda(\ceil{tn})}\le\ord{e^{h(tn+1)}}=\ord{e^{htn}}$ for every $h>0$; hence $o(\Phi^t)\le\inf\EE$. Both inequalities are strict, because $o(\Phi^t)\in\OO$ while $\sup\PP,\inf\EE\notin\OO$. The case $t<0$ follows from Remark~\ref{rem:facts}(iii).
\end{proof}
 
\begin{proof}[Proof of Theorem~\ref{thm:D}]
Given $\varphi$ as in Theorem~\ref{thm:cassaigne}, let $\Lambda_\varphi$ be the minimal two-sided subshift of Section~\ref{sec:symbolic}, $M_\varphi:=M_{\Lambda_\varphi}$, which is a one-dimensional continuum by Proposition~\ref{prop:continuum}, and let $f_\varphi:=\Phi^{t_\varphi}$, $t_\varphi\in(0,1)$, be the minimal homeomorphism given by Corollary~\ref{cor:homeo}, so that $o(f_\varphi)=\ord{P_\varphi(n)}$ with $P_\varphi(n):=p_{\Lambda_\varphi}(\ceil{t_\varphi n})$. Since $\log p_{\Lambda_\varphi}(m)\sim\varphi(m)$ and $\ceil{t_\varphi n}\le2t_\varphi n$ for large $n$, we have
\begin{equation}\label{eq:Pphi}
\tfrac12\varphi(t_\varphi n)\le\log P_\varphi(n)\le2\varphi(2t_\varphi n)\qquad\text{for all large }n.
\end{equation}
Since $\log t/\varphi(t)\to0$, and $\varphi(t)/t\to0$ because $t^{\beta}\varphi'(t)\to0$, we have \[\ord{n^d}\le\ord{p_{\Lambda_\varphi}(n)}\le\ord{e^{hn}}\] for all $d,h>0$, so \[\sup\PP<o(\sigma|_{\Lambda_\varphi})<\inf\EE\] and, by Theorem~\ref{thm:C}(ii), $\sup\PP<o(f_\varphi)<\inf\EE$.
 
(a) Let $\varphi(t)=t^\alpha$ and $f_\alpha:=f_\varphi$. By \eqref{eq:Pphi}, the claimed bounds hold with $c_\alpha=\frac{\ln2}{2}t_\varphi^\alpha$ and $C_\alpha=2\ln2\,(2t_\varphi)^\alpha$, since finitely many $n$ do not affect orders. The inequalities $\sup\PP<o(f_\alpha)<\inf\EE$ were proved above. If $\alpha<\alpha'$, then $e^{C_\alpha n^\alpha}/e^{c_{\alpha'}n^{\alpha'}}\to0$, and Remark~\ref{rem:facts}(ii) gives
\[
o(f_\alpha)\le\ord{e^{C_\alpha n^\alpha}}<\ord{e^{c_{\alpha'}n^{\alpha'}}}\le o(f_{\alpha'}).
\]
 
(b) Let $f:=f_{\varphi_1}$, with $\varphi_1$ as in \eqref{eq:exotic}, and $g:=f_{\varphi_2}$, with $\varphi_2(t)=t^{1/2}$; both have intermediate generalized entropy, as shown above. Write $\theta(t):=\frac12+\frac14\cos(\ln\ln(t+10))$, so that $\varphi_1(t)=(t+10)^{\theta(t)}$. Since $\ln\ln(t+10)$ tends to infinity while its increments on $[t,2t]$ tend to zero, there are arbitrarily large $n$ with $\theta\ge\frac{7}{10}$ on $[t_{\varphi_1}n,2t_{\varphi_1}n]$, and arbitrarily large $n$ with $\theta\le\frac{3}{10}$ on that interval. For the former, \eqref{eq:Pphi} gives $\log P_{\varphi_1}(n)\ge\frac12(t_{\varphi_1}n)^{7/10}$, while $\log P_{\varphi_2}(n)\le2(2n)^{1/2}$; hence $P_{\varphi_1}/P_{\varphi_2}$ is unbounded. For the latter, $\log P_{\varphi_1}(n)\le2(2n+10)^{3/10}$, while $\log P_{\varphi_2}(n)\ge\frac12(t_{\varphi_2}n)^{1/2}$; hence $\liminf P_{\varphi_1}/P_{\varphi_2}=0$. Therefore $o(f)$ and $o(g)$ are incomparable.
 
\end{proof}
 
\section{Final remarks and questions}\label{sec:questions}
 
\subsection{Complexities in $\cB$.} Theorem~\ref{thm:A} realizes the orders $\ord{n\,a(n)}$, $a\in\cB$, a proper subset of the orders in $\cB$ above $\ord{n}$: there are $c\in\cB$ with $c(n)\ge n$ such that $c(n)/n$ is not $\asymp$-equivalent to any nondecreasing sequence. Is every order in $\cB$ above $\ord{n}$ the complexity class of a minimal subshift? Which orders below $\sup\PP$ but outside $\cB$ are generalized entropies of minimal homeomorphisms of continua?

\subsection{Intermediate orders.}
Cassaigne's construction gives $\log p_\Lambda(n)\sim\varphi(n)$, which does not determine the class $\ord{p_\Lambda(n)}$: an error $o(\varphi(n))$ in the logarithm can yield an unbounded multiplicative factor after exponentiation. For $0<\alpha<1$, is the exact order $\ord{e^{n^\alpha}}$ realized by a minimal homeomorphism of a one-dimensional continuum? More generally, which orders strictly between $\sup\PP$ and $\inf\EE$ occur?

\subsection{Low complexity on continua.}
Every nonzero time map of the suspension of an infinite subshift has generalized entropy at least $\ord{n}$, by the Morse--Hedlund theorem and Proposition~\ref{prop:flow}. Minimal circle rotations, on the other hand, have generalized entropy $\ord{1}$. Which orders strictly between $\ord{1}$ and $\ord{n}$ occur for minimal homeomorphisms of one-dimensional continua?

For the particular spaces $M_a$, even the order $\ord{1}$ is excluded for minimal homeomorphisms. Indeed, a minimal isometry would make $M_a$ homogeneous: the compact closure of its iterates in the isometry group acts transitively. By \cite[Theorem~5.2]{CH}, homogeneity would imply equicontinuity of the holonomy pseudogroup of the suspension foliation. This pseudogroup contains the iterates of the infinite subshift $\Lambda_a$, which are expansive and cannot be equicontinuous. Thus $M_a$ admits no minimal isometry for any compatible metric. A minimal homeomorphism with generalized entropy $\ord{1}$ would be equicontinuous \cite[Theorem~3]{CP}, and hence an isometry for a compatible metric, giving the same contradiction. Does every minimal homeomorphism $g$ of $M_a$ satisfy $o(g)\ge o(f_a)$?

\subsection{Variable roof functions.} If $p_\Lambda\in\cB$, is $o(\Phi_r^t)=\ord{p_\Lambda(n)}$ for every continuous roof $r$ and $t\ne0$ (Remark~\ref{rem:roof})?
 
\bibliographystyle{plain}
\bibliography{references}

\end{document}